\documentclass{amsart}
\usepackage[utf8]{inputenc}

\usepackage[demo]{graphicx}
\usepackage{wallpaper}
\usepackage{wrapfig,booktabs}

\usepackage{fancyhdr}
\usepackage{lettrine}

\usepackage{mdframed}

\usepackage{cancel}

\usepackage[margin=1.45in]{geometry}

\usepackage{adjustbox}

\usepackage[nointegrals]{wasysym}

\usepackage{tikz}

\usepackage{amsmath}
\usepackage{mathtools}
\usepackage{amsthm}
\usepackage{amssymb}
\usepackage{tcolorbox}
\tcbuselibrary{theorems}
\usepackage{tikz-cd}
\usepackage{xparse}

\newcommand{\newthought}[1]{\medskip\noindent\fbox{{#1}}}

\newcommand{\abs}[1]{| #1 |}

\newcommand{\cat}[1]{\mathcal{#1}}

\newcommand{\gen}[1]{\left\langle #1 \right\rangle}

\usepackage[colorlinks=true]{hyperref}

\usepackage[all]{hypcap} % When using hyperref this will allow you to jump to images, rather than the captions below them

\usepackage[nameinlink]{cleveref} % Better references to lemmas etc. The nameinlink makes all of 'Theorem 1' clickable, not just the number at the end

\usepackage{tabularx}

\usepackage[backend=biber]{biblatex}
\title[On Diamond-Free Subgroup Lattices]{On Diamond-Free Subgroup Lattices}

\author{Matt Alexander}
\date{\today}

\address{University of Regina\\
3737 Wascana Parkway\\
Regina, Saskatchewan, S4S 0A2\\
Canada}
\email{mpa097@uregina.ca}

\begin{document}

%=================================
% Colourblind-friendly colours (?)
%=================================

% Need more testing to see if actually colourblind-friendly, and for what kinds of colourblindness

\definecolor{lemmacolour}{RGB}{0, 114, 178}

\definecolor{propcolour}{rgb}{0.6, 0.0,0.0}

\definecolor{theoremcolour}{rgb}{1.0, 0.84, 0.0}

\definecolor{examplecolour}{RGB}{59,191,143}

\definecolor{definitioncolour}{rgb}{0.75, 0.75, 0.75}

%=================================
% ENVIRONMENTS
%=================================

\newenvironment{proof-lemma}
    {\begin{mdframed}[linewidth=3, linecolor=lemmacolour,topline=false,rightline=false,bottomline=false]
\begin{proof}
}{
\end{proof}
\end{mdframed}}

\newenvironment{proof-prop}
    {\begin{mdframed}[linewidth=3, linecolor=propcolour,topline=false,rightline=false,bottomline=false]
\begin{proof}
}{
\end{proof}
\end{mdframed}}

\newenvironment{proof-theorem}
    {\begin{mdframed}[linewidth=3, linecolor=theoremcolour,topline=false,rightline=false,bottomline=false]
\begin{proof}
}{
\end{proof}
\end{mdframed}}

\newenvironment{defframe}
    {\begin{mdframed}[linewidth=3, linecolor=definitioncolour,topline=false,rightline=false,bottomline=false]
}{
\end{mdframed}}

\newenvironment{exampleframe}
    {\begin{mdframed}[linewidth=3, linecolor=examplecolour,topline=false,rightline=false,bottomline=false]
}{
\end{mdframed}}

%=================================================================================

\newtcbtheorem[number within=section]{theorem}{\color{black}Theorem}%
{colback=theoremcolour!5,colframe=theoremcolour,fonttitle=\bfseries,
label type=theorem
}{thm}

\crefname{theorem}{theorem}{theorems}

\newtcbtheorem[number within=section, use counter from=theorem]{definition}{Definition}%
{colback=definitioncolour!5,colframe=definitioncolour!50!black,fonttitle=\bfseries,
label type=definition}{def}

\crefname{definition}{definition}{definitions}

\newtcbtheorem[number within=section, use counter from=theorem]{lemma}{Lemma}%
{colback=lemmacolour!5,colframe=lemmacolour!95!black,fonttitle=\bfseries, 
label type=lemma}{lem}

\crefname{lemma}{lemma}{lemmas}

\newtcbtheorem[number within=section, use counter from=theorem]{prop}{Proposition}%
{colback=propcolour!5,colframe=propcolour,fonttitle=\bfseries, 
label type=proposition}{propin}
% Use propin or something like it (I think there's a conflict with the use of "prop", which makes clever ref not like it)

\crefname{proposition}{proposition}{propositions}

\newtcbtheorem[number within=section, use counter from=theorem]{example}{Example}%
{colback=examplecolour!5,colframe=examplecolour!90!black,fonttitle=\bfseries, 
label type=example}{ex}

\crefname{example}{example}{examples}

%============================
%============================

\begin{abstract}

In this paper we introduce a particular lattice of subgroups called a ``cyclic-diamond'' and show that every finite non-cyclic group contains a cyclic-diamond as a sublattice of its lattice of subgroups. Turning to the infinite case, we show that an infinite abelian group does not contain a cyclic-diamond in its subgroup lattice if and only if all of its finitely generated subgroups are cyclic or isomorphic to $\mathbb{Z} \times \mathbb{Z}_{2^N}$ for some $N$.

\end{abstract}

\keywords{Subgroup lattice, distributive lattice, locally cyclic group}
\subjclass[2020]{20E15}

\maketitle

\tableofcontents

%=================================
\section{Introduction}
%=================================

It is known that the subgroup lattice of a group is intimately tied to the structure of the group itself. A number of group-theoretic properties can be reframed in terms of the subgroup lattice. For example, supersolvable groups are equivalently those whose subgroup lattices satisfy the Jordan--Dedekind chain condition (all maximal chains of subgroups have the same length) \cite{iwasawa1941endlichen}.

In \cite{ore1938structures}, Ore proved that the subgroup lattice of a finite group $G$ is distributive if and only if $G$ is cyclic. This was later extended to the fact that an arbitrary group $G$ has a distributive subgroup lattice if and only if $G$ is locally cyclic (i.e. all finitely generated subgroups are cyclic). This latter result can be found, for example, in \cite[Theorem~1.2.3]{schmidt1994subgroup}.

It is a known fact that there are only 2 possible obstructions to a lattice being distributive: a lattice is distributive if and only if it does not contain a copy of $M_3$ (the diamond) or $N_5$ (the pentagon) as a sublattice (see for instance \cite[Chapter~1, Theorem~3.6]{sankappanavar1981course}).

\begin{figure}
    \centering  
\begin{tikzcd}[every arrow/.append style={dash}]
&\bullet \ar[dl] \ar[dr] \ar[d] \\
\bullet \ar[dr] & \bullet \ar[d] & \bullet \ar[dl] \\
& \bullet
\end{tikzcd}
\quad
\begin{tikzcd}[every arrow/.append style={dash}]
&\bullet \ar[dl] \ar[ddr] \\
\bullet \ar[d] & &  \\
\bullet \ar[dr] & & \bullet \ar[dl]\\
& \bullet
\end{tikzcd}
    \caption{The diamond $M_3$ and the pentagon $N_5$.}
    \label{fig:diamond-pentagon}
\end{figure}
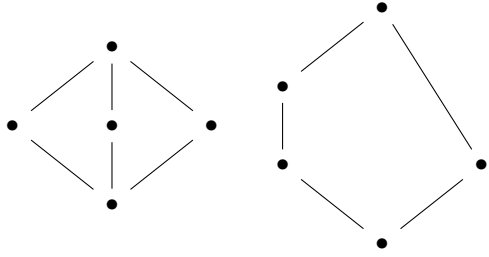

Lattices that lack $N_5$ are called \textbf{modular}, and have been well-studied in the group context: the lattice of normal subgroups of any group is modular. In particular, every abelian group has a modular subgroup lattice. It follows that every abelian non-locally-cyclic group must contain a diamond. In this paper we refine this result by looking at the specific form of diamonds that appear in subgroup lattices.

\newthought{Notation: } Let $\gen{x_1, \ldots, x_n}$ denote the subgroup generated by elements $x_1, \ldots, x_n$.

\begin{theorem}{{\color{black}Main Result}}{}
Every finite non-cyclic group $G$ contains a diamond of the form

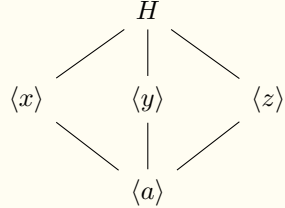
\begin{center}
\begin{tikzcd}[every arrow/.append style={dash}]
& H \ar[dl] \ar[dr] \ar[d] \\
\gen{x} \ar[dr] & \gen{y} \ar[d] & \gen{z} \ar[dl] \\
& \gen{a}
\end{tikzcd}
\end{center}
for some subgroup $H \leq G$,
where $\gen{a}$ is a maximal subgroup of both $\gen{x}$ and $\gen{y}$.
\end{theorem}

\newthought{Terminology: } We will refer to diamonds of the above type as \textbf{cyclic-diamonds} (note that $H$ may or may not be cyclic). Groups which do not contain such diamonds in their subgroup lattices will be called \textbf{cyclic-diamond free}, and in the finite case they are equivalently cyclic groups. However the infinite case is more subtle, as we discuss in \Cref{sec:inf}. Note that every cyclic-diamond is contained in a finitely generated subgroup, as $H$ is generated by any two of the elements $\{x,y,z\}$.

%-----------------------------------------------------------

\subsection*{Acknowledgments} The author would like to thank Martin Frankland for fruitful discussions and feedback and Allen Herman for his helpful discussions. This work was supported by an NSERC Postgraduate Scholarship.

%-----------------------------------------------------------

%=================================================
\section{Finite Group Case}\label{sec:fin}
%=================================================

\begin{figure}
    \centering
    \begin{tikzcd}[every arrow/.append style={dash}]
&Q_8 \ar[dl] \ar[dr] \ar[d] \\
\gen{i} \ar[dr] & \gen{j} \ar[d] & \gen{k} \ar[dl] \\
& \gen{-1} \ar[d]\\
& e
\end{tikzcd}
\quad
\begin{tikzcd}[every arrow/.append style={dash}]
&\mathbb{Z}^2_p \ar[dl] \ar[dr] \ar[d] \ar[drr] \\
\gen{(1,0)} \ar[dr] & \gen{(0,1)} \ar[d] & \gen{(1,1)} \ar[dl] & \ldots \ar[dll]\\
& e
\end{tikzcd}

    \caption{The subgroup lattices of $Q_8$ and $\mathbb{Z}_p^2$.}
    \label{fig:q8-zp2-lattice}
\end{figure}
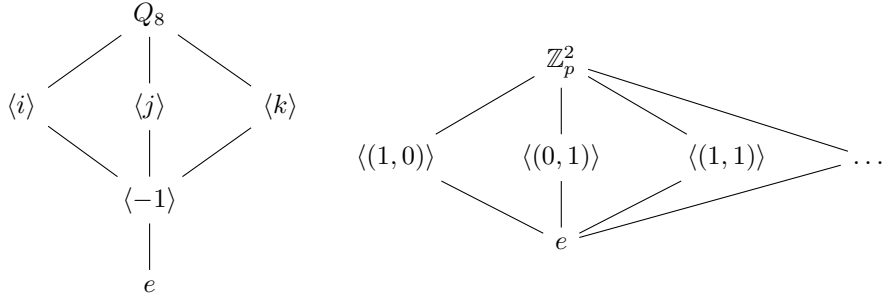

\newthought{Miller--Moreno: } In \cite{miller1903non} it was shown that every non-cyclic group $G$ whose subgroups are all cyclic is isomorphic to one of the following: 

\begin{itemize}
    \item $Q_8$ (the quaternions),

    \item $\mathbb{Z}_p^2$ for some prime $p$,

    \item The semidirect product $\mathbb{Z}_q \rtimes \mathbb{Z}_{p^a}$ with $q \equiv 1$ (mod $p$) (of which there is only one possibility up to isomorphism for each pair of suitable primes $p, q$).
\end{itemize}
Looking at the subgroup lattices of $Q_8$ and $\mathbb{Z}_p^2$ shown in \Cref{fig:q8-zp2-lattice}, we can see that they contain cyclic-diamonds. However, it is not immediately clear what happens in the case of $\mathbb{Z}_q \rtimes \mathbb{Z}_{p^a}$. Below we demonstrate that diamonds of the above form occur within groups of this type as~well. As every finite group contains a minimal non-cyclic subgroup, we will conclude that such diamonds occur within any finite non-cyclic group.

\begin{lemma}{Non-Abelian Prime Generated Case}{}
Let $G$ be non-abelian and generated by its prime order elements. Then $G$ contains a cyclic-diamond of the form

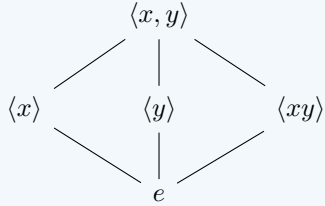
\begin{center}
\begin{tikzcd}[every arrow/.append style={dash}]
&\gen{x,y} \ar[dl] \ar[dr] \ar[d] \\
\gen{x} \ar[dr] & \gen{y} \ar[d] & \gen{xy} \ar[dl] \\
& e
\end{tikzcd}
\end{center}
for every pair of non-commuting prime order elements $x, y$.

\end{lemma}
\begin{proof-lemma}
Note that powers of prime order elements also have prime order (or are the identity). Thus as $G$ is generated by its prime order elements, if all prime order elements commute, then $G$ would be abelian. So there must be two prime order elements: $x,y$ with $\abs{x} = p$ and $\abs{y} = q$ such that  $xy \ne yx$. In particular, this implies that $\gen{x} \ne \gen{y}$ (otherwise $x$ and $y$ would commute). Consider the lattice

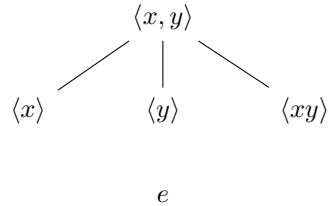
\begin{center}
\begin{tikzcd}[every arrow/.append style={dash}]
&\gen{x,y} \ar[dl] \ar[dr] \ar[d] \\
\gen{x}  & \gen{y}  & \gen{xy} \\
& e
\end{tikzcd}
\end{center}
It's clear that any two of $\{ x, y, xy \}$ generate $\gen{x,y}$, so the upper part of the diamond holds. Since $x$ and $y$ have prime order, they have no non-trivial proper subgroups. Thus $\gen{x} \wedge \gen{y} = e$. 

\newthought{Consider $\gen{x} \wedge \gen{xy}$: } As $x$ has prime order, the only subgroups of $\gen{x}$ are itself and the identity. If $\gen{x} \wedge \gen{xy} = \gen{x}$, we would have $x \in \gen{xy}$, so $x = (xy)^k$ for some $k$. Then $y^{-1} = y^{-1} x^{-1} x = y^{-1}  x^{-1} (xy)^k = (xy)^{k-1}$, demonstrating $y^{-1} \in \gen{xy}$. Thus $x,y \in \gen{xy}$. As $\gen{xy}$ is cyclic, this implies that $x$ and $y$ commute, which does not hold by assumption. Thus 
\begin{equation*}
\gen{x} \wedge \gen{xy} = e
\end{equation*}
By a symmetric argument, $\gen{y} \wedge \gen{xy} = e$, and so we find that $G$ contains a diamond of the claimed form.
\end{proof-lemma}

\phantomsection\label{fix:pa-q}

\begin{prop}{$p^a q$}{pa-q}
Let $G$ be a non-abelian group of order $p^a q$ that is minimal non-cyclic and not prime generated. Then $G$ contains a cyclic-diamond of the form

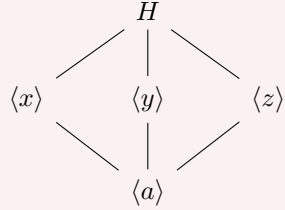
\begin{center}
\begin{tikzcd}[every arrow/.append style={dash}]
& H \ar[dl] \ar[dr] \ar[d] \\
\gen{x} \ar[dr] & \gen{y} \ar[d] & \gen{z} \ar[dl] \\
& \gen{a}
\end{tikzcd}
\end{center}

where $\gen{a}$ is maximal in $\gen{x}$ and $\gen{y}$.
\end{prop}
\begin{proof-prop}
\vphantom{a}

\newthought{Unique subgroup of order $p$: } If $G$ contains two distinct prime subgroups $\gen{x} \ne \gen{y}$ with $\abs{x} = \abs{y} = p$ then $\gen{x,y}$ would not be cyclic (since cyclic groups contain at most one subgroup of any given order). Since $G$ is minimal non-cyclic, this would imply that $G = \gen{x,y}$, but $G$ is not prime generated by assumption. Thus $G$ must have a unique subgroup of order $p$.

\newthought{Larger unique subgroups: } Let $\cat{A} = \{ i \mid 1 \leq i \leq a-1, G$ has a unique subgroup of order $p^i \}$. $\cat{A}$ is non-empty and thus contains a largest element $N < a$. Note that every $i \leq N$ is contained in $\cat{A}$: if $P$ is a subgroup of order $p^i$, it is contained in some Sylow $p$-subgroup, which in this case is cyclic and thus must contain unique subgroups of order $p^i$ and $p^N$. It follows that every subgroup of order $p^i$ with $i \leq N$ must be contained in the unique subgroup of order $p^N$ (which is cyclic), and thus there must be a unique subgroup of order $p^i$. Now let $P_i$ be the unique subgroup of order $p^i$ for $1 \leq i \leq N$. Now $G$ must contain two distinct subgroups of order $p^{N+1}$, call these $\gen{x} \ne \gen{y}$ (they must be cyclic as $G$ is minimal non-cylic, and these are clearly not equal to $G$ as they contain no element of order $q$). Now consider the following lattice:

% Worth mentioning maybe, is the fact that if $i \leq N$ then we must have a unique subgroup of order i. This is because every p-power subgroup is contained in some sylow p group (which in this case is cyclic). So if $i \leq N$, we have some subgroup of order $p^i$, but then it's contained in some sylow subgroup, which has a unique subgroup of each possible order (since it's cyclic). So that Sylow subgroup would have a unique subgroup of order $p^i$ and a unique subgroup of order $p^N$. Thus it must be that that order $p^i$ subgroup is contained in the unique subgroup of order $p^N$, which we know is unique. So we have a unique subgroup of order $p^i$ for all $i \leq N$.

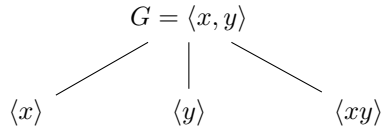
\begin{center}
\begin{tikzcd}[every arrow/.append style={dash}]
& G = \gen{x,y} \ar[dl] \ar[dr] \ar[d] \\
\gen{x} & \gen{y}  & \gen{xy}
\end{tikzcd}
\end{center}

\noindent As $\gen{x}$ and $\gen{y}$ both contain a subgroup of order $p^N$, they must both contain $P_N$. As $P_N$ is maximal in both $\gen{x}$ and $\gen{y}$, if it were not equal to their intersection, then we would have $\gen{x} \wedge \gen{y}$ equal to both $\gen{x}$ and $\gen{y}$, and in particular, $\gen{x} = \gen{y}$, which does~not hold by assumption. Thus
\begin{equation}
\gen{x} \wedge \gen{y} = P_N.
\end{equation}
We next consider the other intersections. The order of $xy$ is of the form $\abs{xy} = p^i q^j$, with $i \leq a$ and $j \leq 1$.

\newthought{We just need $P_N, \gen{xy} \leq H$: } If $H$ is a proper subgroup of $G$ that contains both $\gen{xy}$ and $P_N$, then it must be the case that $\gen{x} \wedge H = \gen{y} \wedge H = P_N$. Since $\gen{x} \wedge H$ is the largest subgroup contained in both $\gen{x}$ and $H$, and $P_N$ is maximal in $\gen{x}$, if $\gen{x} \wedge H \ne P_N$, we would have $\gen{x} \wedge H = \gen{x}$. But as $xy \in H$, this would also imply that $y \in H$. So $\gen{x,y} = G \leq H$, and thus $H = G$. But by assumption, $H$ is a proper subgroup. Thus $\gen{x} \wedge H = P_N$, and by a similar argument $\gen{y} \wedge H = P_N$. Below, we will find such a subgroup $H$.

\newthought{If $\abs{xy} = p^i$: }

\begin{itemize}

\item \fbox{If $i \leq N$: } Then $\gen{xy}$ is the unique subgroup of order $p^i$, and thus $\gen{xy} \leq \gen{x} \wedge \gen{y}$. So $xy = x^k$, implying $y = x^{k-1}$ and $y \in \gen{x}$. Similarly, $x \in \gen{y}$, so $\gen{x} = \gen{y}$, which does not hold by assumption.

\item \fbox{If $i > N$: } Then $\gen{xy}$ contains the unique subgroup of order $p^N$: $P_N \leq \gen{xy}$. As $\gen{xy}$ is a proper subgroup of $G$ (since $G$ isn't cyclic), it follows from an argument above that $P_N = \gen{x} \wedge \gen{xy} = \gen{y} \wedge \gen{xy}$. Thus in this case $G$ contains a diamond of the form

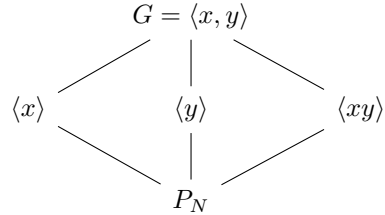
\begin{center}
\begin{tikzcd}[every arrow/.append style={dash}]
& G = \gen{x,y} \ar[dl] \ar[dr] \ar[d] \\
\gen{x} \ar[dr] & \gen{y} \ar[d] & \gen{xy} \ar[dl] \\
& P_N
\end{tikzcd}
\end{center}

\end{itemize}

\newthought{If $\abs{xy} = p^t q$: } From Miller--Moreno, $G$ contains a unique subgroup of order $q$ \cite{miller1903non}, say $Q$. As $q$ divides the order of $xy$, it follows that $Q \leq \gen{xy}$.

\begin{itemize}

\item \fbox{If $t \leq N$: } There is a unique subgroup, $P_t$ of order $p^t$, which must be contained in $\gen{xy}$. Thus $P_t,Q \leq \gen{xy}$, so $P_t \vee Q \leq \gen{xy}$. In fact, $\gen{xy} = P_t \vee Q$ (any subgroup which contains both $P_t$ and $Q$ must have order at least $p^t q = \abs{xy}$, and thus the smallest of these must be $\gen{xy}$).

\newthought{Subgroup of order $p^{a-1}q$: } From \cite{miller1903non}, we know that $G$ contains a cyclic subgroup of order $p^{a-1} q$, say $C_{p^{a-1} q}$. Since $t \leq N \leq a-1$, it follows that $C_{p^{a-1} q}$ contains subgroups of order $p^N$, $p^t$ and $q$, and thus must contain $P_N, P_t$ and $Q$. Thus 
\begin{equation}
\gen{xy} = P_t \vee Q \leq C_{p^{a-1} q}
\end{equation}
We can now show that $G$ contains a diamond of the form

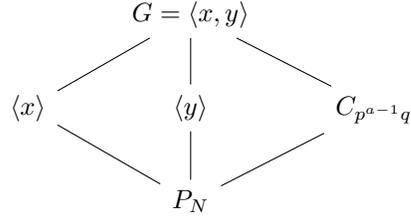
\begin{center}
\begin{tikzcd}[every arrow/.append style={dash}]
& G = \gen{x,y} \ar[dl] \ar[dr] \ar[d] \\
\gen{x} \ar[dr] & \gen{y} \ar[d] & C_{p^{a-1} q} \ar[dl] \\
& P_N
\end{tikzcd}
\end{center}
The lower part of the diamond follows from an argument above, since $P_N, \gen{xy} \leq C_{p^{a-1} q}$. For the upper part, note that since $\gen{xy} \leq C_{p^{a-1} q}$, we have $\gen{x} \vee \gen{xy} = G \leq \gen{x} \vee C_{p^{a-1} q}$, so $G = \gen{x} \vee C_{p^{a-1} q}$, and similarly $G = \gen{y} \vee C_{p^{a-1} q}$.

\item \fbox{If $t > N$: } As $p^N$ divides the order of $xy$, we have $P_N \leq \gen{xy}$, and thus we can apply the same argument as before to show that we have a diamond of the form

\begin{center}
\begin{tikzcd}[every arrow/.append style={dash}]
& G = \gen{x,y} \ar[dl] \ar[dr] \ar[d] \\
\gen{x} \ar[dr] & \gen{y} \ar[d] & \gen{xy} \ar[dl] \\
& P_N
\end{tikzcd}
\end{center}

\end{itemize}

\end{proof-prop}

\phantomsection\label{fix:fin-abelian-case}

\begin{lemma}{Finite Abelian Case}{fin-abelian-case}
Let $G$ be finite and abelian. Then $G$ is cyclic-diamond free if and only if $G$ is cyclic.
\end{lemma}
\begin{proof-lemma}

\newthought{If cyclic: } From Ore \cite{ore1938structures}, finite cyclic groups are equivalently those with distributive subgroup lattices. In particular, such groups are free of all diamonds (including cyclic-diamonds).

\newthought{If cyclic-diamond free: } If $G$ contained two distinct subgroups of the same prime order, $p$, then it would contain a copy of $\mathbb{Z}_p^2$ (whose lattice of subgroups contains a cyclic-diamond), appearing as the join of those subgroups. By the Fundamental Theorem of Finitely Generated Abelian Groups, it follows that $G$ is of the form
\begin{equation}
G \cong \mathbb{Z}_{p_1^{a_1}} \times \ldots \times \mathbb{Z}_{p_n^{a_n}} \cong \mathbb{Z}_{p_1^{a_1} \ldots p_n^{a_n}}
\end{equation}
for distinct primes $p_i$, and is thus cyclic.
\end{proof-lemma}

\newthought{Collecting our results } we have the following:

\phantomsection\label{fix:finite-cyc-diam-free}

\begin{theorem}{{\color{black}Main Result}}{finite-cyc-diam-free}
A finite group is cyclic if and only if its subgroup lattice does not contain a diamond of the form

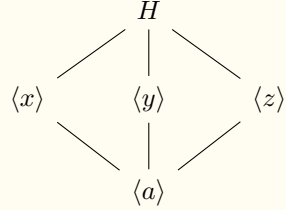
\begin{center}
\begin{tikzcd}[every arrow/.append style={dash}]
& H \ar[dl] \ar[dr] \ar[d] \\
\gen{x} \ar[dr] & \gen{y} \ar[d] & \gen{z} \ar[dl] \\
& \gen{a}
\end{tikzcd}
\end{center}

where $\gen{a}$ is a maximal subgroup of both $\gen{x}$ and $\gen{y}$.
\end{theorem}
\begin{proof-theorem}

\newthought{If cyclic: } By Ore \cite{ore1938structures}, the subgroup lattice is distributive and thus contains no diamond of any type.

\newthought{If not cyclic: } Then as $G$ is finite, it contains a minimal non-cyclic subgroup, $H$. From \cite{miller1903non}, $H$ must be isomorphic to $\mathbb{Z}_p^2$, $Q_8$, or the non-trivial semidirect product $\mathbb{Z}_q \rtimes \mathbb{Z}_{p}^a$. In the first two cases, simple inspection of the subgroup lattice of $H$ shows that it contains a cyclic-diamond. A cyclic-diamond exists in the third case by \hyperref[fix:pa-q]{\Cref{propin:pa-q}}.
% \Cref{propin:pa-q}.
%
\end{proof-theorem}

\begin{prop}{}{}
Let $G$ be a finite group. The following are equivalent:
\begin{enumerate}
\item The subgroup lattice of $G$ is distributive.
\item The subgroup lattice of $G$ is cyclic-diamond free.
\item $G$ is cyclic.
\end{enumerate}
\end{prop}
\begin{proof-prop}
\vphantom{a}
\newthought{$(1) \Rightarrow (2)$: } This is a known lattice-theoretic result: a lattice is distributive if and only if it is diamond-free (in particular, cyclic-diamond free) and modular.

\newthought{$(2) \Rightarrow (3)$: } From \hyperref[fix:finite-cyc-diam-free]{\Cref{thm:finite-cyc-diam-free}} the only finite cyclic-diamond free groups are cyclic.
% \Cref{thm:finite-cyc-diam-free} 

\newthought{$(3) \Leftrightarrow (1)$: } This was shown by Ore \cite{ore1938structures}.
\end{proof-prop}

%------------------------------------------------------
\subsection{Organizing Finitely Generated Groups}
%------------------------------------------------------

\vphantom{a}

In \Cref{tab:groups} we provide examples of finitely generated groups, organized by the structure relevant to the proofs above. That is, such groups are organized based on whether they are minimal non-cylic, whether they are prime generated, and whether there is some prime $p$ such that the group contains a unique subgroup of order $p$. Examples of finite groups are provided unless no such examples exist.

\begin{table}[t]

\begin{center}

\begin{tabularx}{\textwidth}{llll>{\raggedright\arraybackslash}X}
\toprule
\multicolumn{2}{l}{\hphantom{aaaaaa}Group}\\
\vphantom{a}
Abelian & Non-Abelian & Minimal Non-cyclic? & Prime Generated? &  Unique Prime Subgroup?\\
\midrule
None & $S_3$ & Yes & Yes & Yes \\
$\mathbb{Z}_p \times \mathbb{Z}_p$ & Tarski Monsters & Yes & Yes & No \\
None &  $Q_8$ & Yes & No & Yes\\
None & None & Yes & No & No\\
$\mathbb{Z}_p^2 \times \mathbb{Z}_q$ & $S_3 \times \mathbb{Z}_2$  & No & Yes & Yes\\
$\mathbb{Z}_p^2 \times \mathbb{Z}_q^2$ & $A_4$ & No & Yes & No\\
$\mathbb{Z}_{p^2} \times \mathbb{Z}_{pq}$ & $Q_8 \times \mathbb{Z}_3$ & No & No & Yes\\
$\mathbb{Z}_{p^2} \times \mathbb{Z}_p$ & $Q_8 \times \mathbb{Z}_2$ & No & No & No
\\
\bottomrule

\end{tabularx}

\end{center}
    
    \caption{Structure of finitely generated non-cyclic groups.}
    \label{tab:groups}

\end{table}

\newthought{Abelian Minimal Non-cyclic: } From \cite{miller1903non} we know that the only finite abelian minimal non-cyclic groups are of the form $\mathbb{Z}_p \times \mathbb{Z}_p$. No other finitely generated abelian minimal non-cyclic groups can be found, as if $G \cong \mathbb{Z}^r \times \mathbb{Z}_{p_1^{a_1}} \times \ldots \times \mathbb{Z}_{p_n^{a_n}}$ and $r \geq 1$, $G$ will contain non-cyclic proper subgroups isomorphic $\mathbb{Z} \times \mathbb{Z}_n$. Thus only groups of the form $\mathbb{Z}_p \times \mathbb{Z}_p$ are found in the first four rows of the abelian column of \Cref{tab:groups}.

\newthought{Non-Abelian Minimal Non-cyclic: } From \cite{miller1903non}, for every finite non-abelian minimal non-cyclic group $G$ there is some prime $p$ such that $G$ contains a unique subgroup of order $p$. Thus there are no finite non-abelian examples of groups of type Yes-Yes-No or Yes-No-No in \Cref{tab:groups}. However, there are \textit{finitely generated} non-abelian groups of type Yes-Yes-No: Tarski Monsters of type $p$, which were originally shown to exist in \cite{ol1980infinite}, are infinite groups in which every non-trivial subgroup has order $p$. In particular, such groups are minimal non-cylic, prime generated (being generated by any two elements contained in different subgroups), and do not have a unique subgroup of a given prime order. On the other hand, there are no groups of type Yes-No-No: any such group would contain at least two distinct subgroups of order $p$ for some prime $p$, say $\gen{x}$ and $\gen{y}$.  It would follow that $G = \gen{x,y}$ (being minimal non-cyclic), and so $G$ would also be prime generated.

%====================================================
\section{Infinite Abelian Groups}\label{sec:inf}
%====================================================

 % \Cref{lem:fin-abelian-case} 

We can try to apply a similar argument to \hyperref[fix:fin-abelian-case]{\Cref{lem:fin-abelian-case}} in the case that $G$ is an arbitrary abelian group. However, once $G$ becomes infinite, the maximality condition on cyclic-diamonds becomes too strong of a requirement. If we remove that condition, we will be able to classify abelian groups in terms of their (generalized) cyclic-diamonds.

\newthought{Terminology: } In what follows, a \textbf{generalized-cyclic-diamond} will refer to a sublattice of the form

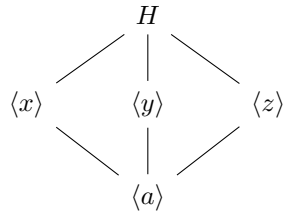
\begin{center}
\begin{tikzcd}[every arrow/.append style={dash}]
& H \ar[dl] \ar[dr] \ar[d] \\
\gen{x} \ar[dr] & \gen{y} \ar[d] & \gen{z} \ar[dl] \\
& \gen{a}
\end{tikzcd}
\end{center}
\textit{without} any conditions on the maximality of the bottom subgroup.

\newthought{Abelian Groups: } Any finitely generated abelian group can be expressed (up to isomorphism) as $\mathbb{Z}^r \times \mathbb{Z}_{p_1^{a_1}} \times \ldots \times \mathbb{Z}_{p_n^{a_n}}$, where $r$ is the rank of the group. The argument in \hyperref[fix:fin-abelian-case]{\Cref{lem:fin-abelian-case}} can be applied to show that the primes $\{ p_i \mid 1 \leq i \leq n\}$, must be distinct. Additionally, if the rank of the group is at least two, its subgroup lattice would contain a generalized-cyclic-diamond of the form

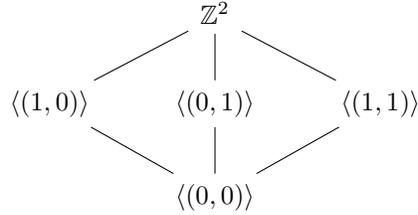
\begin{center}
\begin{tikzcd}[every arrow/.append style={dash}]
&\mathbb{Z}^2 \ar[dl] \ar[dr] \ar[d] \\
\gen{(1,0)} \ar[dr] & \gen{(0,1)} \ar[d] & \gen{(1,1)} \ar[dl] \\
& \gen{(0, 0)}
\end{tikzcd}
\end{center}
(note that the bottom subgroup here is not maximal in any of the middle subgroups). Thus abelian generalized-cyclic-diamond free groups must be such that all of their finitely generated subgroups are cyclic or of the form $\mathbb{Z} \times \mathbb{Z}_n$. Groups of the latter form are abelian and not locally-cyclic (as they are finitely generated and not cyclic), and thus must contain diamonds. However, it is not obvious that such groups must contain \textit{cyclic-diamonds}. In what follows we will show that $\mathbb{Z} \times \mathbb{Z}_n$ is generalized-cyclic-diamond free if and only if $n = 2^N$ for some $N$.

%---------------------------------------------------------------------------------------------------
\subsection{Groups of the Form \texorpdfstring{$G \cong \mathbb{Z} \times \mathbb{Z}_n$}{TEXT} }
%---------------------------------------------------------------------------------------------------

\vphantom{a}

\newthought{$p \mid n$: } If $n$ is divisible by some prime greater than $2$, we can find a cyclic-diamond in the subgroup lattice of $G$.

% \newthought{Non-locally cyclic: } $G$ contains countably infinite many copies of itself generated as $\gen{(n,0)} \vee \gen{(0,1)}$ for all $n \geq 2$. Thus $G$ is far from being locally cyclic. However, as $G$ is abelian, its subgroup lattice cannot contain a pentagon \cite{ore1937structures}, and so must contain a diamond somewhere in its subgroup lattice.

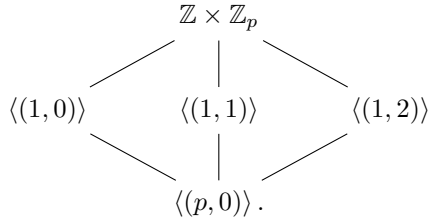
\phantomsection\label{fix:z-zn-diamond}

\begin{example}{Generalized-Cyclic-Diamond When $p \mid n$, $p > 2$}{z-zn-diamond}
Consider $G = \mathbb{Z} \times \mathbb{Z}_n$ where $p \mid n$, for some prime $p > 2$. It follows that $G$ contains a subgroup $H \cong \mathbb{Z} \times \mathbb{Z}_{p}$. Such a subgroup will contain a diamond of the form

\end{example}
\begin{exampleframe}

\begin{center}
\begin{tikzcd}[every arrow/.append style={dash}]
& \mathbb{Z} \times \mathbb{Z}_p \ar[dl] \ar[dr] \ar[d] \\
\gen{(1,0)} \ar[dr] & \gen{(1,1)} \ar[d] & \gen{(1,2)} \ar[dl] \\
& \gen{(p,0)}.
\end{tikzcd}
\end{center}

\end{exampleframe}

\newthought{$\mathbb{Z} \times \mathbb{Z}_{2^N}$: } We reduce to the case where $n$ is a power of $2$. Note that we cannot treat this case by applying the above construction, as $\gen{(1,2)} = \gen{(1,0)}$ in $\mathbb{Z} \times \mathbb{Z}_2$. Surprisingly, we will see that $\mathbb{Z} \times \mathbb{Z}_{2^N}$ is actually generalized-cyclic-diamond free!

\newthought{Subgroups: } Note that the proper subgroups of $\mathbb{Z} \times \mathbb{Z}_{2^N}$ are (up to isomorphism) $\mathbb{Z}$, $\mathbb{Z}_{2^k}$, and $\mathbb{Z} \times \mathbb{Z}_{2^k}$, with $1 \leq k \leq N$. As $\mathbb{Z}$ and $\mathbb{Z}_{2^k}$ are cyclic, they contain no diamonds (including generalized-cyclic-diamonds). Thus if it can be shown that diamonds of the form
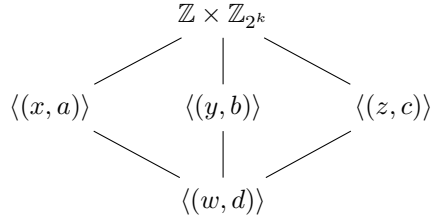
\begin{center}
\begin{tikzcd}[every arrow/.append style={dash}]
& \mathbb{Z} \times \mathbb{Z}_{2^k} \ar[dl] \ar[dr] \ar[d] \\
\gen{(x,a)} \ar[dr] & \gen{(y,b)} \ar[d] & \gen{(z,c)} \ar[dl] \\
& \gen{(w,d)}
\end{tikzcd}
\end{center}
where $0 \leq a,b,c,d < 2^k$ cannot exist, it will follow that all groups of the form $\mathbb{Z} \times \mathbb{Z}_{2^N}$ are generalized-cyclic-diamond free. For the sake of space, we will give diamonds of the above shape a special name.

\phantomsection\label{fix:model-cyclic-diamond}

\begin{definition}{Generalized-Cyclic-Diamond of Model Type}{model-cyclic-diamond}

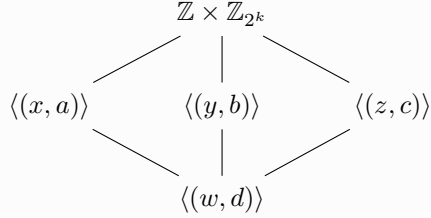
\begin{center}
\begin{tikzcd}[every arrow/.append style={dash}]
& \mathbb{Z} \times \mathbb{Z}_{2^k} \ar[dl] \ar[dr] \ar[d] \\
\gen{(x,a)} \ar[dr] & \gen{(y,b)} \ar[d] & \gen{(z,c)} \ar[dl] \\
& \gen{(w,d)}
\end{tikzcd}
\end{center}
where $0 \leq a,b,c,d < 2^k$.
\end{definition}
\begin{defframe}
\label{defin-frame}
\newthought{Bezout: } If a diamond of the above shape were to exist, the element $(1,0)$ could be written as a linear combination of any two of $\{ (x,a), (y,b), (z,c) \}$, and in particular, the integer $1$ could be written as a linear combination of any two of $\{x,y,z \}$. From Bezout's theorem it follows that $\{ x,y,z \}$ are pairwise coprime.

\end{defframe}

\newthought{Notation: } Denote the \textbf{$2$-adic valuation} of an integer $x$ by $\nu_2(x)$.

\phantomsection\label{fix:abc-even}

\begin{lemma}{}{abc-even}
In a diamond of the form of \hyperref[fix:model-cyclic-diamond]{\Cref{def:model-cyclic-diamond}} at least two of $\{a, b, c \}$ are odd.
\end{lemma}
\begin{proof-lemma}
Without loss of generality say $a$ and $b$ are both even. Then any linear combination $\lambda a + \mu b$ is even, and will thus remain even mod $2^N$. It follows that $(0,1)$ is not an element of $\gen{(x,a)} \vee \gen{(y,b)} = \mathbb{Z} \times \mathbb{Z}_{2^N}$, which cannot be. Thus at most one element of $\{a, b, c \}$ is even.
\end{proof-lemma}

\phantomsection\label{fix:2-adic-match}

\begin{lemma}{}{2-adic-match}
In a diamond of the form of \hyperref[fix:model-cyclic-diamond]{\Cref{def:model-cyclic-diamond}} where $d \ne 0$, we have $\nu_2(i) = \nu_2(j)$ for each ordered pair $(i,j) \in \{ (x,a), (y,b), (z,c), (w,d) \}$.

\end{lemma}
\begin{proof-lemma}
As $\gen{(w,d)} \leq  \gen{(x,a)}$, we have that $w = \alpha x$ and $d \equiv \alpha a$ (mod $2^N$) for some $\alpha \in \mathbb{Z}$. Thus $\nu_2(w) = \nu_2(\alpha) + \nu_2(x)$. Note that we can write $\alpha a = d + 2^N k$ for some $k \in \mathbb{Z}$. As $1 \leq d < 2^N$, $0 \leq \nu_2(d) \leq N-1$, and thus we can write
\begin{equation*}
\alpha a = 2^{\nu_2(d)} (\dfrac{d}{2^{\nu_2(d)}} + 2^{N - \nu_2(d)} k)
\end{equation*}
Now $d/2^{\nu_2(d)}$ is odd, and $2^{N - \nu_2(d)} k$ is even (since $\nu_2(d) < N$). Thus the term inside of the brackets is odd. It follows that $\nu_2(\alpha a) = \nu_2(d)$. A similar result holds for $b$ and $c$.

\newthought{$\nu_2(w) = \nu_2(d)$: } As $\{x,y,z\}$  are coprime, two of these elements must be odd. Also from \hyperref[fix:abc-even]{\Cref{lem:abc-even}}, we have that two of $\{ a,b,c \}$ are odd. It follows that at least one of $\{ (x,a), (y,b), (z,c) \}$ must have both of its components odd. Say $y$ and $b$ are both odd.  Write $w = \beta y$ and $d \equiv \beta b$ (mod $2^N$) for some $\beta$. We have $\nu_2(\beta b) = \nu_2(d)$, but as $b$ is odd, $\nu_2(\beta b) = \nu_2(\beta)$. Now as $y$ is odd, we have $\nu_2(w) = \nu_2(\beta y) = \nu_2(\beta) = \nu_2(d)$.

\newthought{$\nu_2(x) = \nu_2(a)$: } As $y$ and $b$ are both odd, we have $\nu_2(y) = \nu_2(b) = 1$. Now as we have $w = \alpha x$, $\nu_2(w) = \nu_2(\alpha x) = \nu_2(\alpha) + \nu_2(x)$. But $\nu_2(w) = \nu_2(d) = \nu_2(\alpha a) = \nu_2(\alpha) + \nu_2(a)$. Thus $\nu_2(\alpha) + \nu_2(x) = \nu_2(\alpha) + \nu_2(a)$ and so $\nu_2(x) = \nu_2(a)$. A similar argument shows $\nu_2(z) = \nu_2(c)$.
\end{proof-lemma}

\phantomsection\label{fix:both-comp-odd}

\begin{lemma}{}{both-comp-odd}
In a diamond of the form of \hyperref[fix:model-cyclic-diamond]{\Cref{def:model-cyclic-diamond}} at least two of $\{ (x,a), (y,b), (z,c) \}$ have both of their components odd.

\end{lemma}
\begin{proof-lemma}

As $\{ x,y,z \}$ are coprime, at least two of them are odd. Say $y$ and $z$ are odd. From \hyperref[fix:abc-even]{\Cref{lem:abc-even}} at least two of $\{ a,b,c \}$ must be odd as well. It follows that at least one of $\{ (x,a), (y,b), (z,c) \}$ must have both of its components odd. Without loss of generality, say $y$ and $b$ are odd. If $c$ is odd as well, we can consider the pairs $\{ (y,b), (z,c) \}$. Thus we reduce to considering the case where $c$ is even (and thus $a$ must be odd). If $x$ is odd, we have the pairs $\{ (x,a), (y,b) \}$. Thus the only remaining case to consider is where $\{x,c \}$ are even and $\{a,y,b,z \}$ are odd.

\newthought{If $d \ne 0$: } Then from \hyperref[fix:2-adic-match]{\Cref{lem:2-adic-match}} we have that $\nu_2(b) = \nu_2(y)$ and $\nu_2(c) = \nu_2(z)$ and thus $\{ (y,b), (z,c) \}$ have both of their components odd.

\newthought{If $d = 0$: } Consider $\gen{(x,a)} \wedge \gen{(y,b)} = \gen{(w,d)} = \gen{(w,0)}$. Note that $\gen{(2^N x, 0)} \leq \gen{(x,a)}$, as $t (2^N x, 0) = 2^N t (x,a)$. Thus $\gen{(2^N x,0)} \wedge \gen{(2^N y,0)} \leq \gen{(x,a)} \wedge \gen{(y,b)} = \gen{(w,0)}$. On the other hand, $w = \alpha x = \beta y$ and $\alpha a \equiv \beta b \equiv 0$ (mod $2^N$) for some $\alpha, \beta$. As $a,b$ are odd, $\alpha a \equiv \beta b \equiv 0$ (mod $2^N$) implies that $\alpha = 2^N k$ and $\beta = 2^N m$ for some $k,m$. Thus $w = 2^N k x = 2^N m y$, and so $(w,0) \in \gen{(2^N x, 0)} \wedge \gen{(2^N y, 0)}$. Thus
\begin{equation}
\gen{(w,0)} = \gen{(2^N x,0)} \wedge \gen{(2^N y,0)}
\end{equation}
As everything is living in a copy of $\mathbb{Z}$, it follows that $\gen{2^N x} \wedge \gen{2^N y} = \gen{w}$, and thus $w = \pm \operatorname{lcm}(2^N x, 2^N y) = \pm 2^N \operatorname{lcm}(x,y)$. As $x$ and $y$ are coprime, we have 
\begin{equation}
w = \pm 2^N xy
\end{equation}
Now $(w,0) \in \gen{z,c}$ implies $w = \pm 2^N xy = \gamma z$ for some $\gamma$. As $\{x,y,z \}$ are coprime, this implies $z = \pm 2^r$ for some $r \geq 0$. Since $z$ is odd, $z = \pm 1$. Now $(2^N y, 0) = \pm 2^N y (\pm 1,c) = 2^N (y,b)$ so $(2^N y, 0) \in \gen{(z,c)} \wedge \gen{(y,b)} = \gen{(w,d)} = \gen{(2^N xy, 0)}$. So $2^N y = \lambda 2^N xy$ for some $\lambda$. It follows that $x = \pm 1$, but $x$ is even. Thus this case cannot occur.
\end{proof-lemma}

\phantomsection\label{fix:z-zn-cd-free}

\begin{theorem}{{\color{black}$\mathbb{Z} \times \mathbb{Z}_{2^N}$ is Generalized-Cyclic-Diamond Free}}{z-zn-cd-free}
The subgroup lattice of $\mathbb{Z} \times \mathbb{Z}_{2^N}$ does not contain a generalized-cyclic-diamond.

\end{theorem}
\begin{proof-theorem}

From \hyperref[fix:both-comp-odd]{\Cref{lem:both-comp-odd}} we have that at least two of $\{ (x,a), (y,b), (z,c) \}$ have both of their components odd. Without loss of generality, say $\{ y,b,z,c \}$ are all odd. Consider $z(y,b) = (yz,zb)$ and $y(z,c) = (yz, yc)$.

\newthought{Multiply by $2^{N-1}$: } Let $g$ and $h$ be odd integers. Then their difference is even: $g-h = 2m$. Thus $2^{N-1} g - 2^{N-1} h = 2^N m \equiv 0$ (mod $2^N$). So $2^{N-1}g \equiv 2^{N-1} h$ (mod $2^N$). And as $g$ and $h$ are odd, and thus contain no factors of $2$, we have $2^{N-1}g , 2^{N-1} h \not\equiv 0$ (mod~$2^N$). Now $2^{N-1}z(y,b) = (2^{N-1} yz, 2^{N-1} zb)$ and $2^{N-1}y(z,c) = (2^{N-1} yz, 2^{N-1} yc)$. As $\{y,z,b,c \}$ are all odd, $yc$ and $zb$ are odd, and thus $2^{N-1} zb \equiv 2^{N-1} yc$ (mod $2^N$). It follows that $(2^{N-1} yz, 2^{N-1} zb) = (2^{N-1} yz, 2^{N-1} yc)$ as elements of $\mathbb{Z} \times \mathbb{Z}_{2^N}$. Thus
\begin{equation*}
(2^{N-1} yz, 2^{N-1} zb) \in \gen{(y,b)} \wedge \gen{(z,c)} = \gen{(w,d)}
\end{equation*}
Note that, as $z$ and $b$ are both odd, $2^{N-1} zb \not\equiv 0$ (mod $2^N$). So in particular we have that $d \ne 0$ and thus we can apply the results of \hyperref[fix:2-adic-match]{\Cref{lem:2-adic-match}} to conclude $\nu_2(x) = \nu_2(a)$. 

\newthought{If $x = 1$: } Since $\nu_2(x) = \nu_2(a)$, $a$ must also be odd. If $z = 1$ as well, then writing $(0,1) = \lambda (x,a) + \mu (z,c) = \lambda(1,a) + \mu(1,c)$ we would have $\mu = -\lambda$ and so $\lambda a + \mu c = \lambda (a - c) \equiv 1$ (mod $2^N$). But $a$ and $c$ are both odd, and thus their difference is even, and remains even mod $2^N$. It follows that $\lambda (a - c) \not\equiv 1$ (mod $2^N$). So we must have $z \ne 1$ (similarly $z \ne -1$). So if $x =1$, we can redo the above procedure using $\{ x, a, y, b \}$ instead of $\{ y,b, z, c \}$ to make sure that the element of $\{ x,y,z \}$ that we do not use in our construction is not $\pm 1$. So without loss of generality, we can say $x \ne \pm 1$.

\newthought{$x = 2^k$: } As $\gen{(w,d)} \leq \gen{(x,a)}$, we have that $(2^{N-1} yz, 2^{N-1} zb) \in \gen{(x,a)}$, and thus $2^{N-1} yz = \lambda x$ for some $\lambda$. But $\{x,y,z \}$ are coprime, and so $x$ must be of the form $x = \pm 2^k$ for some $1 \leq k \leq N-1$. In particular, $x$ is even. As $\nu_2(x) = \nu_2(a)$, we have that $a$ is even as well.

\newthought{$\gen{(x,a)} \vee \gen{(y,b)}$: } Consider $(0,1) \in \gen{(x,a)} \vee \gen{(y,b)} = \mathbb{Z} \times \mathbb{Z}_{2^N}$. We must have $(0,1) = \lambda(x,a) + \mu(y,b)$ for some $\lambda, \mu$. So $\lambda x = - \mu y$. As $x$ and $y$ are coprime and $x$ is even, $\mu$ must be even. Thus $\lambda a + \mu b$ is even, and will remain even mod $2^N$. But $\lambda a + \mu b \equiv 1$ (mod $2^N$). Thus we see that a diamond of the type of \hyperref[fix:model-cyclic-diamond]{\Cref{def:model-cyclic-diamond}} cannot be formed, and so $\mathbb{Z} \times \mathbb{Z}_{2^N}$ is generalized-cyclic-diamond free.
\end{proof-theorem}

\begin{example}{Non-Locally-Cylic but Generalized-Cyclic-Diamond Free}{}
The groups $\mathbb{Z} \times \mathbb{Z}_{2^N}$ for $N \in \mathbb{N}$ are non-locally-cyclic abelian groups whose subgroup lattice contains neither a copy of the pentagon nor a diamond of the form 

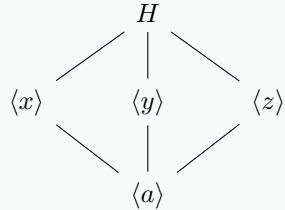
\begin{center}
\begin{tikzcd}[every arrow/.append style={dash}]
& H \ar[dl] \ar[dr] \ar[d] \\
\gen{x} \ar[dr] & \gen{y} \ar[d] & \gen{z} \ar[dl] \\
& \gen{a}
\end{tikzcd}
\end{center}
(even if $\gen{a}$ is non-maximal in $\gen{x}, \gen{y}, \gen{z}$).
\end{example}
\begin{exampleframe}
However, as expected, $\mathbb{Z} \times \mathbb{Z}_{2^N}$ will contain non-cyclic diamonds. For example:

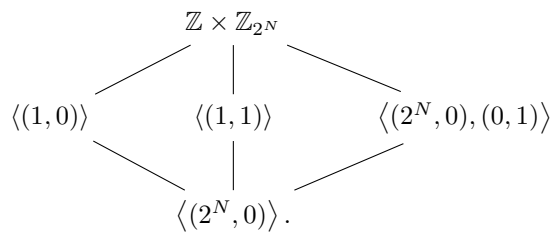
\begin{center}
\begin{tikzcd}[every arrow/.append style={dash}]
& \mathbb{Z} \times \mathbb{Z}_{2^N} \ar[dl] \ar[dr] \ar[d] \\
\gen{(1,0)} \ar[dr] & \gen{(1,1)} \ar[d] & \gen{(2^N,0), (0,1)} \ar[dl] \\
& \gen{(2^N,0)}.
\end{tikzcd}
\end{center}

\end{exampleframe}

\begin{theorem}{{\color{black}Classification of Abelian Generalized-Cyclic-Diamond Free Groups}}{}
Let $G$ be an abelian group. Then $G$ is generalized-cyclic-diamond free if and only if every finitely generated subgroup $H \leq G$ is cyclic or isomorphic to $\mathbb{Z} \times \mathbb{Z}_{2^N}$ for some $N$.
\end{theorem}
\begin{proof-theorem}

\newthought{If generalized-cyclic-diamond free: } As we remarked at the beginning of \Cref{sec:inf}, any cyclic-diamond free finitely generated abelian group must be cyclic or of the form $\mathbb{Z} \times \mathbb{Z}_n$. From \hyperref[fix:z-zn-diamond]{\Cref{ex:z-zn-diamond}} the only possibilities for $n$ are $2^N$ for some~$N$, proving this direction.

\newthought{Other direction: } By \hyperref[fix:z-zn-cd-free]{\Cref{thm:z-zn-cd-free}}, every $\mathbb{Z} \times \mathbb{Z}_{2^N}$ is generalized-cyclic-diamond free. Since every generalized-cyclic-diamond is contained in a finitely generated subgroup, if the only finitely generated subgroups are cyclic or $\mathbb{Z} \times \mathbb{Z}_{2^N}$, $G$ must be generalized-cyclic-diamond free.
\end{proof-theorem}

\begin{center}
   \line(1,0){200}
 \end{center}

%========================================================
%========================================================

\printbibliography

@article{ore1938structures,
    AUTHOR = {Ore, Oystein},
     TITLE = {Structures and group theory. {II}},
   JOURNAL = {Duke Math. J.},
  FJOURNAL = {Duke Mathematical Journal},
    VOLUME = {4},
      YEAR = {1938},
    NUMBER = {2},
     PAGES = {247--269},
      ISSN = {0012-7094,1547-7398},
   MRCLASS = {99-04},
  MRNUMBER = {1546048},
       DOI = {10.1215/S0012-7094-38-00419-3},
       %URL = {https://doi.org/10.1215/S0012-7094-38-00419-3},
}

@book{schmidt1994subgroup,
    AUTHOR = {Schmidt, Roland},
     TITLE = {Subgroup lattices of groups},
    SERIES = {De Gruyter Expositions in Mathematics},
    VOLUME = {14},
 PUBLISHER = {Walter de Gruyter \& Co., Berlin},
      YEAR = {1994},
     PAGES = {xvi+572},
      ISBN = {3-11-011213-2},
   MRCLASS = {20E15 (20-02 20D30)},
  MRNUMBER = {1292462},
MRREVIEWER = {Francesco\ de Giovanni},
       DOI = {10.1515/9783110868647},
       %URL = {https://doi.org/10.1515/9783110868647},
}

@book{sankappanavar1981course,
    AUTHOR = {Burris, Stanley and Sankappanavar, H. P.},
     TITLE = {A course in universal algebra.},
    SERIES = {},
 PUBLISHER = {Springer-Verlag, New York-Berlin,, },
      YEAR = {1981},
     PAGES = {xvi+276},
      ISBN = {0-387-90578-2},
   MRCLASS = {08-01},
  MRNUMBER = {648287},
MRREVIEWER = {R.\ S.\ Pierce},
}

@article{iwasawa1941endlichen,
    AUTHOR = {Iwasawa, Kenkiti},
     TITLE = {\"{U}ber die endlichen {G}ruppen und die {V}erb\"{a}nde ihrer {U}ntergruppen.},
   JOURNAL = {J. Fac. Sci. Imp. Univ. Tokyo Sect. I.},
  FJOURNAL = {},
      YEAR = {1941},
     PAGES = {171--199},
   MRCLASS = {20.0X},
  MRNUMBER = {5721},
MRREVIEWER = {R.\ Baer},
}

@article{miller1903non,
    AUTHOR = {Miller, G. A. and Moreno, H. C.},
     TITLE = {Non-abelian groups in which every subgroup is abelian},
   JOURNAL = {Trans. Amer. Math. Soc.},
  FJOURNAL = {Transactions of the American Mathematical Society},
    VOLUME = {4},
      YEAR = {1903},
    NUMBER = {4},
     PAGES = {398--404},
      ISSN = {0002-9947,1088-6850},
   MRCLASS = {20D99},
  MRNUMBER = {1500650},
       DOI = {10.2307/1986409},
       %URL = {https://doi.org/10.2307/1986409},
}

@article{ol1980infinite,
    AUTHOR = {Ol'\v{s}anski\u{i}, A. Ju.},
     TITLE = {An infinite group with subgroups of prime orders.},
   JOURNAL = {Izv. Akad. Nauk SSSR Ser. Mat.},
  FJOURNAL = {Izvestiya Akademii Nauk SSSR. Seriya Matematicheskaya},
      YEAR = {1980},
    NUMBER = {no. 2,},
     PAGES = {309--321, 479},
      ISSN = {0373-2436},
   MRCLASS = {20E07},
  MRNUMBER = {571100},
MRREVIEWER = {N.\ S.\ Chernikov},
}

\end{document}